\documentclass[11pt,oneside,notitlepage]{article}

\usepackage[english]{babel}
\usepackage{amsmath}
\usepackage{amssymb}
\usepackage{amsthm}
\usepackage{mathtools}
\usepackage{xspace}

\usepackage{fullpage}

\newtheorem{theorem}{Theorem}

\makeatletter
\let\@fnsymbol\@arabic
\makeatother

\def\BB{\mathbb B}

\def\RR{\mathbb R}

\def\cH{\mathcal H}

\def\cH{\mathcal H}

\begin{document}

\title{
  On Quadratization of Pseudo-Boolean Functions\footnote{
    This paper appeared in the special session on Boolean and pseudo-Boolean functions
    of the 2012 edition of the International Symposium on Artificial Intelligence and 
    Mathematics (ISAIM) held in Fort Lauderdale, FL, USA, on January 9--11, 2012.
  }
}

\author{
  Endre Boros\thanks{MSIS and RUTCOR, Rutgers University, NJ, USA. \emph{E-mail:} endre.boros@rutgers.edu}
  \qquad
  Aritanan Gruber\thanks{RUTCOR, Rutgers University, NJ, USA. \emph{E-mail:} aritanan.gruber@rutgers.edu}
}

\maketitle

\begin{abstract}
  We survey current term-wise techniques for quadratizing high-degree pseudo-Boolean functions and introduce
  a new one, which allows multiple splits of terms. We also introduce the first aggregative approach, which 
  splits a collection of terms based on their common parts.
\end{abstract}

\section{Introduction}

Set functions, i.e., real mappings form the family of subsets of a finite set to the reals are known and widely used in discrete mathematics for almost a century, and in particular in the last 50 years. If we replace a finite set with its characteristic vector, then the same set function can be interpreted as a mapping from the set of binary vectors to the reals. Such mappings are called pseudo-Boolean functions and were introduced in the works of Peter L. Hammer in the 1960s, see the seminal book \cite{HR68}. Pseudo-Boolean functions are different from set functions, only in the sense that their algebraic representation, a multilinear polynomial expression, is usually assumed to be available as an input representation:
\begin{equation}\label{e-pBf}
f(x_1,x_2,\ldots, x_n) ~=~ \sum_{S\subseteq V} a_S\prod_{j\in S}x_j,
\end{equation}
where $V=\{1,2,...,n\}$ is the set of variable indices, and where we assume that $x_j\in \BB=\{0,1\}$ for all $j\in V$. Due to this assumption we have $x_j^2=x_j$, and hence any polynomial expression in such binary variables is indeed equivalent with a multilinear one. In fact it was shown in \cite{HR68} that any set function has a unique multilinear polynomial representation. The \emph{degree} of $f$ is defined as the size of the largest set $S$ with a nonzero coefficient in the above (unique) multilinear polynomial representation of $f$: $\deg(f)~=~ \min \{|S|\mid S\subseteq V,~ a_S\neq 0\}$. Clearly, the degree of a constant function is zero. We say that $f$ is a \emph{quadratic} (resp. \emph{linear}) pseudo-Boolean function, if $\deg(f)\leq 2$ (resp. $\deg(f)\leq 1$).

The problem of minimizing a pseudo-Boolean function (over the set of binary vectors) appears to be the common form of numerous optimization problems, including the well-known MAX-SAT and MAX-CUT problems, and have applications in areas ranging from physics through chip design to computer vision; see e.g., the surveys \cite{BH02,KZ04,RB09}.

Some of these applications lead to the minimization of a quadratic pseudo-Boolean function, and hence such quadratic binary optimization problems received ample attention in the past decades. The survey \cite{BH02} describes a large set of computational tools for such problems. One of the most frequently used technique is based on roof-duality \cite{HHS84}, and aims at finding in polynomial time a simpler form of the given quadratic minimization problem, by fixing some of the variables at their provably optimum value (persistency) and decomposing the residual problem into variable disjoint smaller subproblems, see \cite{BH89,BHST08}. The method in fact was found very effective in computer vision problems, where frequently it can fix up to 80-90\% of the variables at their provably optimum value. This algorithm was recoded by computer vision experts and a very efficient implementation, called QPBO, is freely downloadable, see \cite{RKLS07}.

In many applications of pseudo-Boolean optimization the objective function \eqref{e-pBf} is a higher degree multilinear polynomial. For such problems there are substantially fewer effective techniques available. In particular, there is no analogue to the persistencies (fixing variables at their provably optimum value) provided by roof-duality for the quadratic case. On the other hand, more and more applications would demand efficient methods for the minimization of such higher degree pseudo-Boolean functions; see e.g., \cite{Han79,BR07,RB09}. This increased interest, in particular in the computer vision community, lead to a systematic study of methods to reduce a higher degree minimization problem to a quadratic one; see e.g., \cite{FD05,Ish09,Ish11,KZ04,Ros75,RKFJ09,RRLT11,ZCJ09,ZCJ09b,ZJ10}.

\bigskip

In this paper first we recall known ``quadratization'' techniques from the literature. Next, we provide several new techniques for quadratization, analyze their effectiveness, and recall recent computer vision applications demonstrating their usefulness \cite{FGBZ11}.

\bigskip

\section{Basic Model and Literature Review}

Given a pseudo-Boolean function $f:\BB^n\rightarrow\RR$ as in \eqref{e-pBf}, where $\RR$ denotes the set of reals, the following minimization problem
\begin{equation}\label{e-main}
\min_{x\in\BB^n} f(x)
\end{equation}
is the common form of numerous combinatorial optimization problems. To reduce the above problem to a quadratic minimization problem, we are looking for a quadratic pseudo-Boolean function $g(x,w)$, where $w\in\BB^m$ is a set of ``new'' variables, such that the equality
\begin{equation}\label{e-model}
f(x) ~=~ \min_{w\in\BB^m} g(x,w)
\end{equation}
holds for all $x\in \BB^n$. We shall call such a $g$ the \emph{quadratization} of $f$. The major objective in the problem of quadratizing a given pseudo-Boolean function $f$ is to find such a quadratic function $g$ that satisfies \eqref{e-model}. Secondary objectives are the minimization of the number of new variables $m$, and the ``submodularity" of $g$.

Given two binary vectors $x,y\in \BB^V$, we define their \emph{disjunction} and \emph{conjunction} respectively by $(x\vee y)_j=x_j\vee y_j$ and $(x\wedge y)_j=x_j\wedge y_j$ for all indices $j\in V$. Then, we call a pseudo-Boolean function $f(x)$ \emph{submodular}, if
\[
f(x\vee y) + f(x \wedge y) \leq f(x) + f(y)
\]
holds for any two vectors $x,y\in\BB^V$. Submodular functions play an important role in optimization, since problem \eqref{e-main} which is NP-hard in general, is known to be polynomially solvable if $f$ is submodular \cite{GLS81,IFF00,Sch00}. Let us add that if $f$ is a quadratic pseudo-Boolean function, then it is submodular if and only if all quadratic terms have nonpositive coefficients. This property leads to a very simple, network flow based minimization algorithm \cite{Ham65}. In fact the QPBO implementation returns automatically a minimizing solution for submodular inputs. A similarly efficient characterization of submodularity for cubic pseudo-Boolean functions was also given in \cite{BM85}. Recognition of submodularity for pseudo-Boolean functions of degree $4$ or higher was shown to be NP-hard \cite{Cra89,GS89}.

Let us note that if a pseudo-Boolean function $f$ given in \eqref{e-pBf} can be quadratized \eqref{e-model} by a submodular quadratic function $g$, then $f$ itself must be submodular. It is however not obvious which submodular pseudo-Boolean functions can be quadratized by submodular quadratic functions. As we shall see later, it is easy to show that any cubic submodular function can be quadratized keeping submodularity. However, recent results in  \cite{ZCJ09,ZCJ09b} show that certain degree 4 submodular pseudo-Boolean functions cannot be represented as in \eqref{e-model} by a submodular quadratic function. This is so, even though we do not limit the number of new variables in \eqref{e-model}.

In this paper we study techniques to find functions $g$ satisfying \eqref{e-model} for a given pseudo-Boolean function $f$  with the purpose that such a functions should be ``small'' and ``easy'' to minimize. More precisely, we shall compare techniques by evaluating the number of new variables, the number of terms, and the number of positive quadratic terms, which is a vague measure of non-submodularity.

\bigskip

\subsection{Literature Review}

Let us first recall the quadratization method suggested by \cite{Ros75}, based on the idea of traditional penalty functions. This method replaces a product $xy$ of two binary variables by a new binary variable $w$ (and hence decreases the degree by one of all terms involving both $x$ and $y$), and adds to $f$ a quadratic penalty function $p(x,y,w)$ such that
\begin{equation}\label{e-rosenberg}
p(x,y,w)~ \left\{\begin{array}{ll}=~0&\text{ if } w=xy,\\ \geq ~1&\text{ otherwise.}\end{array}\right.
\end{equation}
Since $f$ is multilinear, we can write it as $f=xyA + B$, where $A$ is a multilinear polynomial not involving $x$ and $y$, and where $B$ is a multilinear polynomial not involving the product $xy$. Assume now that $p$ is a quadratic function satisfying \eqref{e-rosenberg}, and $M$ is a positive real with $M > \max |A|$, where the maximization is taken over all binary assignments of the variables of $A$.
Then, the function $\widetilde{f}~=~ zA + B +Mp$ on $n+1$ variables has the same minima as $f$. \cite{Ros75} showed that
\[
p(x,y,w) ~=~ xy -2xw-2yw+3w
\]
is a quadratic function satisfying \eqref{e-rosenberg}. The above idea then can be applied recursively, until the resulting function $\widetilde{f}$ becomes quadratic. It is easy to see that this is a polynomial transformation, e.g., one never needs more than $O(n^{2\log d})$ new variables, where $d$ is the degree of $f$. The drawback of this approach is that the resulting quadratic function has many ``large'' coefficients, due to the recursive application of the ``big $M$'' substitution. It also introduces many positive quadratic terms, even if the input $f$ is a nice submodular function. These two effects make the minimization of the resulting $\widetilde{f}$ a hard problem, even in approximative sense.

Of course, it would be a simpler approach to replace the product $xy\cdots z$ of several variables by a new variable $w$ and enforce the equality $xy\cdots z=w$ by a quadratic penalty function. It is easy to see that this is not possible with more than two variables.

Let us also note that finding with this approach a quadratization with the minimum number of variables is itself an NP-hard problem. To see this let us consider the cubic pseudo-Boolean function
\[
f(x_0,x_1,...,x_n) ~=~ \sum_{(i,j)\in E} x_0x_ix_j
\]
where $E$ is the edge set of a graph $G$ on vertex set $V=\{1,2,...,n\}$. It is easy to verify that for any quadratization of $f$, there is one with no more new variables, in which we substitute by new variables only products of the form $x_0x_i$, $i\in C$ for a subset $C\subseteq V$ and the ``optimal'' quadratization corresponds to a minimum size vertex cover $C$ of $G$.

\bigskip

A simple quadratization of negative monomials was introduced recently by \cite{KZ04} for degree $3$ monomials, and by \cite{FD05} for arbitrary degree monomials.
\begin{equation}\label{e-KZFD}
-x_1x_2\cdots x_d ~=~ \min_{w\in\BB} ~w\left((d-1) -\sum_{j=1}^d x_j\right).
\end{equation}
Remarkably, all quadratic terms have negative coefficients in this transformation. In particular, if $f$ involves only negative higher degree terms, then the application of \eqref{e-KZFD} yields a submodular quadratization of it. It is also a nice feature of this approach that it does not introduces ``large'' coefficients. Let us also remark that for every minimizing assignment of $f$ there is a corresponding minimum of the quadratized form $g$ obtained by the repeated applications of \eqref{e-KZFD} satisfying $w=x_1x_2\cdots x_d$. Thus, the above transformation can also be viewed as a substitution of a higher degree product.

Let us note that this transformation can be extended easily for positive monomials, as well. For this note first that the equality in \eqref{e-KZFD} is based only on the fact that the symbols ``$x_i$'' stand for a binary value. Thus, introducing \emph{negated literals} $\overline{x}_i=1-x_i$ we can write
\begin{equation}\label{e-KZFD+}
\begin{array}{c}
x_1x_2\cdots x_d - x_{d-1}x_d ~=~\displaystyle  -\sum_{i=1}^{d-2}\overline{x}_i\prod_{j=i+1}^d x_j\\
~=~ \displaystyle \min_{w\in \BB^{d-2}} \sum_{i=1}^{d-2} w_i(d-i - \overline{x}_i-\sum_{j=i+1}^d x_j)
\end{array}
\end{equation}
Hence \eqref{e-KZFD} implies a quadratization of positive monomials as well. For a degree $d$ term we need $d-2$ new variables and get $d-1$ positive (non-submodular) quadratic terms.

Let us add that if a subset of the variables are negated on the left in \eqref{e-KZFD} then the corresponding quadratic terms will have positive coefficients on the right hand side. In particular, if all variables are negated, then all quadratic terms have positive coefficients. Since the new variable $w$ does not appear elsewhere we can replace it with its negation, not changing the minimization in this way, and get again a submodular quadratization
\begin{equation}\label{e-KZFD+n}
\begin{array}{c}
-\overline{x}_1\overline{x}_2\cdots \overline{x}_d ~=~\displaystyle \min_{w\in\BB} ~w\left((d-1) -\sum_{j=1}^d \overline{x}_j\right)\\
~=~\displaystyle -1+\sum_{j=1}^d x_j +\min_{w'\in \BB} w'\left(1-\sum_{j=1}^d x_j\right)
\end{array}
\end{equation}
where $w'=1-w$. Let us call a pseudo-Boolean function $f$ a \emph{unary negaform} if it can be represented as a negative combination of terms involving either only unnegated variables, or only negated ones. It is easy to show that unary negaforms are submodular, and in fact \eqref{e-KZFD} and \eqref{e-KZFD+n} provides a submodular quadratization for such functions.

Unary negaforms (more precisely, their negations) were considered by \cite{BM85} and they showed that all cubic submodular functions can be represented by unary negaforms. They also provided a network flow model for the minimization of a unary negaform. The above \eqref{e-KZFD}, \eqref{e-KZFD+n} submodular quadratization also leads to a network flow based minimization by the results of \cite{Ham65} and these two network flow models are of very similar size (though they are not identical). Thus, the above observations can be viewed as a new simple proof for the results of \cite{BM85}.

Let us remark finally that higher degree submodular functions cannot typically be represented as unary negaforms. This is implied e.g., by the results of \cite{ZCJ09,ZCJ09b} since we just proved that a unary negaform always has a submodular quadratization.

\bigskip

The trick that \eqref{e-KZFD} can be extended by using negated variables was also observed by \cite{RKFJ09} (they called it type-II transformation). They also noticed that one can apply \eqref{e-KZFD} to a subproduct (of a monomial), under some conditions. In particular, they quadratized separately the negated and unnegated variables in a monomial and derived a new transformation (called type-I):
\begin{equation}\label{e-RKFJ}
\small
-\prod_{j\in S_0}\overline{x}_j\prod_{j\in S_1}x_j = \min_{u,v\in \BB} -uv+u\sum_{j\in S_0}x_j +v\sum_{j\in S_1}\overline{x}_j.
\end{equation}

\bigskip

No matter which variation from above we use to quadratize a degree $d$ positive monomial, we always need at least $d-2$ new variables.
In recent publications \cite{Ish09,Ish11} provided a more compact quadratization for positive monomials, using only about half as many variables as the previous methods. To formulate this result, let us consider the positive term $t(x)=x_1x_2\cdots x_d$ of degree $d$, set $k=\lfloor \frac{d-1}{2}\rfloor$, and consider new binary variables $w=(w_1,w_2,...,w_k)$. Define 
\[
  S_1=\sum_{j=1}^d x_j,\quad S_2=\sum_{1\leq i<j\leq d}x_ix_j,\quad
  A=\sum_{j=1}^k w_j\quad\text{and}\quad B=\sum_{j=1}^k (4j-1)w_j.
\]
 Then the following equalities hold:

\begin{equation}\label{e-I-even}
\prod_{j=1}^d x_j ~=~ S_2 + \min_{w\in\BB^k} B-2A S_1
\end{equation}
if $d=2k+2$, and
\begin{equation}\label{e-I-even}
\prod_{j=1}^d x_j ~=~ S_2 + \min_{w\in\BB^k} B-2A S_1 +w_k\left( S_1-d+1\right)
\end{equation}
if $d=2k+1$.

Let us note that $S_1$ and $S_2$ are symmetric functions of $x$, while $A$ is a symmetric function of $w$. However, $B$ is not a symmetric function of $w$. It is an interesting question on its own if one could find a quadratization of $t(x)$ which is symmetric in both $x$ and $w$, and needs substantially fewer new variables than $d$.

Let us also note that while this method introduces substantially fewer variables than the previous methods, it also introduces $\binom{d}{2}$ positive quadratic terms which makes the resulting quadratization highly non-submodular. Despite of this negative feature, \cite{Ish11} reported very good computational results, in particular when compared to the quadratization of \cite{Ros75}.

\bigskip

Let us conclude this section by pointing out that all of the above methods, except \cite{Ros75} introduce individual new variables for each of the monomials of $f$. In many applications this is a disadvantage, increasing the size and frequently the level of non-submodularity of the resulting  quadratization.

We can also note that the quadratization of negative terms are quite well solved by \eqref{e-KZFD}, since we need only one new variable (per term) and the resulting quadratic form is submodular.

In the sequel we present some new quadratization techniques, and for the above reasons, we focus primarily on positive terms and/or on the issue of using fewer than one per term new variables.

\bigskip

\section{Multiple Splits of Terms}

We show first a generalization of some of the above results. We still focus on a single term, and introduce a general scheme to split this term into several fragments in order to decrease the maximum  degree.

Let $p,q$ be positive integers, and denote by $[q]=\{1,2,...,q\}$ the set of positive integers up to $q$. Assume that $\phi_i:\BB^p\rightarrow \BB$ are Boolean functions for $i\in [q]$ satisfying the following conditions:
\begin{equation}\label{en1}
\small
\begin{array}{c}
\displaystyle \min_{y\in \BB^p} \sum_{i=1}^q \phi_i(y) ~=~ 1, ~~\text{ and }\\*[5mm]
\displaystyle \forall I\subseteq [q], ~I\neq [q], ~~\exists y_I\in\BB^p ~~\text{ s.t. }~~ \sum_{i\in I}\phi_i(y_I) ~=~ 0.
\end{array}
\end{equation}
In other words, the sum of the $\phi$ functions have a positive minimum, but if we leave out any of the summands, the minimum becomes zero. For instance, for $p=2$, $q=3$ the functions $\phi_1=y_1$, $\phi_2=y_2$ and $\phi_3=\overline{y}_1\overline{y}_2$ form such a set.

\begin{theorem}\label{tn1}
Let $\phi_i$ be Boolean functions satisfying condition \eqref{en1}, and $P_i\subseteq [d]$ be subsets for $i\in [q]$ covering $[d]$. Then we have
\begin{equation}\label{en2}
\prod_{j=1}^d x_j ~=~ \min_{y\in\BB^p}\sum_{i=1}^q \phi_i(y)\prod_{j\in P_i}x_j.
\end{equation}
\end{theorem}

\proof
If $\prod_{j\in P_i}x_j=1$ for all $i\in [q]$ then we have
\[
1 ~=~ \min_{y\in \BB^p}\sum_{i=1}^q \phi_i(y) ~=~ 1
\]
by \eqref{en1}.  If there is an index $k\in [q]$ for which $\prod_{j\in P_k}x_j=0$, then by \eqref{en1} there exists a $y^*\in \BB^p$ such that $\phi_{i}(y^*)=0$ for all $i\neq k$,  and consequently we have $\phi_k(y^*)=1$. Thus,
\[
\begin{array}{rl}
0 &\displaystyle \leq \min_{y\in \BB^p}\sum_{i=1}^q \phi_i(y)\prod_{j\in P_i}x_j\\
&\displaystyle \leq \sum_{i=1}^q \phi_i(y^*)\prod_{j\in P_i}x_j
~=~ \prod_{j\in P_k}x_j ~=~ 0
\end{array}
\]
follows, proving the claim. \qed

\noindent \textbf{Remarks and Examples}:

\begin{itemize}
\item $\phi_1=y_1$ and $\phi_2=\overline{y}_1=1-y_1$ provides a 2-split;
\item $\phi_1=y_1$, $\phi_2=y_2$, and $\phi_3=\overline{y}_1\overline{y}_2$ provides a 3-split;
\item any binary tree of depth $p$ with $q$  leaves defines an appropriate system of $\phi_i$ functions, however not all systems correspond to such a tree (see e.g., the above 3-split);
\item $p$ variables can in general provide a $q \leq 2^p$-split transforming a degree $d$ term to $q$ terms of maximum degree $p+\lceil\frac{d}{q}\rceil$.
\item a 2-split combined with \eqref{e-KZFD} yields the following quadratization of a cubic term
\[
\begin{array}{rl}
xyz &=~ \displaystyle \min_{u\in \BB} xu+\overline{u}yz\\
&=~ \displaystyle \min_{u\in \BB} xu +yz -uyz\\
&=~ \displaystyle \min_{u,v\in\BB}xu +yz +(2-y-u-z)v;
\end{array}
\]
\item combining the above with a 2-split yields the following quadratization of a quartic term
$txyz = \displaystyle \min_{u\in \BB} txu+\overline{u}yz
= \displaystyle \min_{u,v,w,s\in \BB} tv + xu + (2-v-x-u)w +yz + (2-u-y-z)s$;
another way doing this is
$txyz = \displaystyle \min_{u\in \BB} tu +xyz -uxyz
= \displaystyle \min_{u,v,w,s\in\BB} tu + xv + yz + (2-v-y-z)w + (3-u-x-y-z)s$.
\end{itemize}

It is interesting that in all of the above attempts to quadratize a positive degree $d$ term, we had to include at least $d-1$ positive quadratic terms. We in fact conjecture that this is necessary.

\bigskip

\section{Splitting of Common Parts}

In this section first we still focus on positive terms, and introduce a quadratization which associates a single new variable with several terms, achieving a simultaneous decrease in their degrees.

\begin{theorem}\label{tn2}
Let $C\subseteq [n]$, $\cH\subseteq 2^{[n]\setminus C}$, and consider a fragment of a pseudo-Boolean function of the form
\[
\phi ~=~ \sum_{H\in \cH}\alpha_H\prod_{j\in H\cup C}x_j,
\]
where $\alpha_H\geq 0$ for all $H\in\cH$. Then we have
\begin{equation}\label{en3}
\small
\phi ~=~ \min_{w\in \BB} \left(\sum_{H\in\cH}\alpha_H\right)\overline{w}\prod_{j\in C}x_j + \sum_{H\in\cH} \alpha_H ~w\prod_{j\in H}x_j.
\end{equation}
\end{theorem}

\proof
We claim that
$w=\prod_{j\in C}x_j$ at a minimum, in which  case the left and right hand sides are identical. To see this claim, observe that we have the inequalities
\[
\phi \leq  \sum_{H\in \cH}\alpha_H\prod_{j\in H}x_j,
\]
and
\[
\phi \leq  \sum_{H\in \cH}\alpha_H\prod_{j\in C}x_j.
\]
Furthermore, these two right hand sides are the values of the right hand side of \eqref{en3} corresponding to $w=1$ and $w=0$, respectively. Thus, $w=\prod_{j\in C}x_j$ indeed achieves a value not larger than any of those.
\qed

\bigskip

We extend the idea of splitting away common parts with a single new variable top negative terms, as well.

\begin{theorem}\label{tn3}
Let $C\subseteq [n]$, $\cH\subseteq 2^{[n]\setminus C}$, and consider a fragment of a pseudo-Boolean function of the form
\[
\phi ~=~ -\sum_{H\in \cH}\alpha_H\prod_{j\in H\cup C}x_j,
\]
where $\alpha_H\geq 0$ for all $H\in\cH$. Then we have
\begin{equation}\label{en3+}
\small
\phi ~=~ \min_{w\in \BB} \sum_{H\in\cH} \alpha_H~ w\left(1-\prod_{j\in C}x_j -\prod_{j\in H}x_j\right).
\end{equation}
\end{theorem}

\proof
We can prove, similarly to the previous proof that
$w=\prod_{j\in C}x_j$ at a minimum. For this let us note that if $w=\prod_{j\in C}x_j$, then the right hand side in \eqref{en3+} is identical with $\phi$, since $\prod_{j\in C}x_j\left( 1-\prod_{j\in C}x_j\right) ~=~0$ for all assignments $x\in \BB^n$. Furthermore, we have the inequalities
\[
\phi ~\leq~ 0,
\]
and
\[
\phi ~\leq~ \sum_{H\in\cH} \alpha_H~ \left(1-\prod_{j\in C}x_j -\prod_{j\in H}x_j\right),
\]
where the right hand side values are the right hand side values of \eqref{en3+} corresponding to  $w=0$ and $w=1$, respectively. Thus, again $w=\prod_{j\in C}x_j$ achieves the smallest possible value.
\qed

\bigskip

\section{Conclusions}

In this paper we proposed new quadratization techniques, possibly decreasing the number of new variables needed, when compared to earlier, term-wise quadratization techniques, without introducing ``large'' coefficients and/or unnecessarily many non-submodular terms.

In fact, using Theorem \ref{tn2} recursively, we can find a quadratization $g(x,w)$ of any function $f$ (even if $f$ is already quadratic) in which we have at most $n-1$ positive quadratic terms, where $n$ is the number of variables of $f$. This shows that the difficulty of minimizing $f$ is not coming from the excessive number of non-submodular terms. In fact if we restrict our input to quadratic pseudo-Boolean functions in $n$ variables and with at most $n-1$ positive terms, the minimization problems remains as hard as general quadratic minimization.

On the positive side a recent publication \cite{FGBZ11} demonstrates that for a large class of binary optimization problems arising from computer vision problems quadratization based primarily on Theorem \ref{tn2} is very effective. When compared to the recent results of \cite{Ish11} we obtained quadratizations with substantially fewer new variables and positive quadratic terms. Applying e.g., the polynomial preprocessing algorithm QPBO we managed to run faster and fix substantially more variables at their optimum values than in \cite{Ish11}.

There are several interesting open ends. The most basic one perhaps is a better and more complete understanding of quadratization techniques. Note that in the above results, we did not plan to restrict the number of new variables. It just happened in each case that the resulting quadratized function involved only polynomially many new variables. Is this always the case? Can we get better/easier to minimize quadratizations in terms of exponentially many new variables?

Which submodular functions have submodular quadratization? How to recognize those? How to find efficiently such a quadratization?

\textbf{Acknowledgements.} 
The authors gratefully acknowledge the partial support by NSF grants IIS 0803444 and by CMMI 0856663.
The second author also gratefully acknowledge the partial support by the joint CAPES (Brazil)/Fulbright (USA) 
fellowship process BEX-2387050/15061676.


\end{document}